\documentclass [12pt]{article}

\usepackage{amsfonts}
\usepackage{amssymb}
\usepackage{epsfig, subfigure, verbatim}

\renewcommand{\emph}[1]{{\it #1}}
\newcounter{segcount}

\newenvironment{segment}
{\refstepcounter{segcount}\vspace{5mm}
\noindent{\bf \thesegcount. }}
{}

\newenvironment{statementnumbered}[3][]
{\refstepcounter{segcount}\vspace{5mm}
\noindent{\bf\thesegcount. #2}#1{\bf.}\ {\sl #3}}
{\nolinebreak[4] \nopagebreak[4] $\hfill \square$}

\newenvironment{statement}[3][]
{\vspace{5mm}\noindent{\bf #2}#1{\bf.} {\sl #3}}
{\nolinebreak[4] \nopagebreak[4] $\hfill \square$}

\newenvironment{statementnoboxnumbered}[3][]
{\refstepcounter{segcount}\vspace{5mm}
\noindent{\bf\thesegcount. #2}#1{\bf.} {\sl #3}}
{}

\newenvironment{statementnobox}[3][]
{\vspace{5mm}\noindent{\bf #2}#1{\bf.} {\sl #3}}
{}

\newenvironment{definitionnumbered}[2][]
{\refstepcounter{segcount}\vspace{5mm}
\noindent{\bf\thesegcount. #2}#1{\bf.}}
{}

\newenvironment{definition}[2][]
{\vspace{5mm}\noindent{\bf #2}#1{\bf.}}
{}

\newenvironment{resultnumbered}[3][]
{\refstepcounter{segcount}\vspace{5mm}
\noindent{\bf\thesegcount. #2}#1{\bf.} {\sl #3}
\vskip5mm\noindent {\bf Proof: }}
{\nolinebreak[4] \nopagebreak[4] $\hfill \square$}

\newenvironment{result}[3][]
{\vspace{5mm}
\noindent{\bf#2}#1{\bf.} {\sl #3}
{\\ \bf Proof: }}
{\nolinebreak[4] \nopagebreak[4] $\hfill \square$}

\newenvironment{risultnumbered}[3][]
{\refstepcounter{segcount}\vspace{5mm}
\noindent{\bf\thesegcount. #2}#1{\bf.} {\sl #3}
{\nopagebreak[4] \noindent \bf Proof: }}
{$\hfill \square$}

\newenvironment{risult}[3][]
{\vspace{5mm}
\noindent{\bf#2}#1{\bf.} {\sl #3}
{\noindent \bf Proof: }}
{\nolinebreak[4] \nopagebreak[4] $\hfill \square$}

\newenvironment{risultnoboxnumbered}[3][]
{\refstepcounter{segcount}\vspace{5mm}
\noindent{\bf\thesegcount. #2}#1{\bf.} {\sl #3}
{\nopagebreak[4] \noindent \bf Proof: }}
{}

\newenvironment{risultnobox}[3][]
{\vspace{5mm}
\noindent{\bf#2}#1{\bf.} {\sl #3}
{\nopagebreak[4] \noindent \bf Proof: }}
{}

\newenvironment{resultnoboxnumbered}[3][]
{\refstepcounter{segcount}\vspace{5mm}
\noindent{\bf\thesegcount. #2}#1{\bf.} {\sl #3}
{\\ \bf Proof: }}
{}

\newenvironment{resultnobox}[3][]
{\vspace{5mm}\noindent{\bf #2}#1{\bf.} {\sl #3}
{\\ \bf Proof: }}
{}

\newcommand{\seg}{\begin{segment}}
\newcommand{\segend}{\end{segment}}

\newcommand{\stmtnum}{\begin{statementnumbered}}
\newcommand{\stmtnumend}{\end{statementnumbered}}

\newcommand{\stmt}{\begin{statement}}
\newcommand{\stmtend}{\end{statement}}

\newcommand{\stmtnoboxnum}{\begin{statementnoboxnumbered}}
\newcommand{\stmtnoboxnumend}{\end{statementnoboxnumbered}}

\newcommand{\stmtnobox}{\begin{statementnobox}}
\newcommand{\stmtnoboxend}{\end{statementnobox}}

\newcommand{\defnnum}{\begin{definitionnumbered}}
\newcommand{\defnnumend}{\end{definitionnumbered}}

\newcommand{\defn}{\begin{definition}}
\newcommand{\defnend}{\end{definition}}

\newcommand{\resnum}{\begin{resultnumbered}}
\newcommand{\resnumend}{\end{resultnumbered}}

\newcommand{\res}{\begin{result}}
\newcommand{\resend}{\end{result}}

\newcommand{\risnum}{\begin{risultnumbered}}
\newcommand{\risnumend}{\end{risultnumbered}}

\newcommand{\ris}{\begin{risult}}
\newcommand{\risend}{\end{risult}}

\newcommand{\risnoboxnum}{\begin{risultnoboxnumbered}}
\newcommand{\risnoboxnumend}{\end{risultnoboxnumbered}}

\newcommand{\risnobox}{\begin{risultnobox}}
\newcommand{\risnoboxend}{\end{risultnobox}}

\newcommand{\resnoboxnum}{\begin{resultnoboxnumbered}}
\newcommand{\resnoboxnumend}{\end{resultnoboxnumbered}}

\newcommand{\resnobox}{\begin{resultnobox}}
\newcommand{\resnoboxend}{\end{resultnobox}}

\title{Dually Vertex-Oblique Graphs}

\author{
Alastair Farrugia
\\ {\it afarrugia@math.uwaterloo.ca}
\\ Dept. of Combinatorics \& Optimization
\\ University of Waterloo, Ontario, Canada, N2L 3G1
}

\begin{document}

\maketitle

\begin{abstract}
A vertex with neighbours of degrees 
$d_1 \geq \cdots \geq d_r$ has {\em vertex type} $(d_1, \ldots, d_r)$. A graph is {\em vertex-oblique} if each vertex has a distinct vertex-type. While no graph can have distinct degrees, Schreyer, Walther and Mel'nikov [Vertex oblique graphs, same proceedings] have constructed infinite classes of {\em super vertex-oblique} graphs, where the degree-types of $G$ are distinct even from the degree types of $\overline{G}$.

$G$ is vertex oblique iff $\overline{G}$ is; but $G$ and $\overline{G}$ cannot be isomorphic, since self-complementary graphs always have non-trivial automorphisms. However, we show by construction that there are {\em dually vertex-oblique graphs} of order $n$, where the vertex-type sequence of $G$ is the same as that of $\overline{G}$; they exist iff $n \equiv 0$ or $1 \pmod 4, n \geq 8$, and for $n \geq 12$ we can require them to be split graphs. 

We also show that a dually vertex-oblique graph and its complement are never the unique pair of graphs that have a particular vertex-type sequence; but there are infinitely many super vertex-oblique graphs whose vertex-type sequence is unique.
\end{abstract}

\section{Introduction and basic results}

Let $G$ be a simple graph on $n$ vertices. A vertex $v$ of degree $r$, with neighbours of degrees%
\footnote{It is conventional in the literature on degree sequences to list degrees in non-increasing order. We follow this convention here, even though we do not prefer it, because we will discuss degree sequences in Section~\ref{sec-deg-seq}; it is of little importance anyway.\label{foot-1}}
$x_1 \geq \cdots \geq x_r$, has {\em vertex type} $t(v) := (x_1, \ldots, x_r)$. $G$ is {\em vertex-oblique} if each vertex has a distinct vertex-type.

The degree of $v$ in $\overline{G}$ (the complement of $G$) is $\overline{r} := n-1-r$. If the degrees of vertices in $G$ are $x_1 \geq \cdots \geq x_r, r, y_1 \leq \cdots \leq y_{\overline{r}}$, then $v$ is non-adjacent to vertices of degrees $y_1 \leq \cdots \leq y_{\overline{r}}$, so its vertex-type in $\overline{G}$ is $\overline{t}(v) = (\overline{y_1}, \ldots, \overline{y_{\overline{r}}})$. Thus $G$ is vertex oblique if and only if $\overline{G}$ is.

While no graph can have distinct degrees, Schreyer et al.~\cite{sch} have constructed infinite classes of vertex-oblique graphs. In fact, their examples are {\em super vertex-oblique}, with the degree-types of $G$ being distinct even from degree types of $\overline{G}$. 

It is natural to ask whether there are any self-complementary vertex-oblique graphs, but this is impossible because a self-complementary graph always has non-trivial automorphisms~\cite{mg, na}, obtained by applying twice an(y) isomorphism that maps the graph to its complement. However, in this article we construct infinitely many {\em dually vertex-oblique graphs}, where the set of vertex-types of $G$ is the same as that of $\overline{G}$.

Many simple results for self-complementary graphs still hold under the weaker assumption that $G$ and $\overline{G}$ have the same set of vertex-types. In particular, dually vertex-oblique graphs of order $n$ can only exist if $n \equiv 0$ or $1 \pmod 4$. The main result of this paper is that they exist for all feasible $n$ at least $8$:

\stmtnum{Theorem}
{Dually vertex-oblique graphs of order $n>1$ exist iff $n$ is congruent to $0$ or $1 \pmod 4$, and $n \geq 8$.
}
\stmtnumend 
\newline


We will make use of the following elementary lemma that is inspired by similar results on self-complementary graphs.

\risnum{Lemma\label{basic-lem}}
{Let $G$ be a graph with the same degree sequence as $\overline{G}$, say $d_1 \geq \cdots \geq d_n$. Then:
\begin{itemize}
\item[A.] $d_i + d_{n-i+1} = n-1$, for $i = 1, \ldots, n$.
\item[B.] $n \equiv 0$ or $1 \pmod 4$.
\end{itemize}
If, moreover, $G$ has the same set of vertex-types as $\overline{G}$, and there are $r_d$ vertices of degree $d$, and $s_{x,y}$ edges joining vertices of degrees $x$ and $y$, then:
\begin{itemize}
\item[C.] $s_{y, y} + s_{\overline{y}, \overline{y}} = \frac{\scriptstyle 1}{\scriptstyle 2}{r_y \choose 2}$, and if $x \not= y$, $s_{x, y} + s_{\overline{x}, \overline{y}} = r_xr_y$; in particular, $s_{d,\overline{d}} = \frac{\scriptstyle 1}{\scriptstyle 2} r_d^{\,2}$ for all $d$, except if $d = \overline{d} = (n-1)/2$ in which case $s_{d,\overline{d}}  = \frac{\scriptstyle 1}{\scriptstyle 2}{r_d \choose 2}$.
\item[D.] $r_d$ is even for all $d$, except for $r_{(n-1)/2} \equiv 1 \pmod 4$.
\end{itemize}
Furthermore, if $G$ is dually vertex-oblique, then:
\begin{itemize}
\item[E.] $r_d < 2d, d \not= (n-1)/2$; in particular, there are no isolated or end-vertices.
\item[F.] There must be at least three different degrees in $G$.
\end{itemize}
}
A. If $r_d$ vertices have degree $d$ in $G$, then $r_d$ vertices have degree $d$ in $\overline{G}$ and, thus, degree $n-1-d$ in $G$. So the degree sequence is symmetric about $\frac{\scriptstyle 1}{\scriptstyle 2}(n-1)$.

\noindent B. The number of edges of $G$ and $\overline{G}$ is the same: $\frac{\scriptstyle 1}{\scriptstyle 2}{n \choose 2} = \frac{\scriptstyle 1}{\scriptstyle 4}n(n-1)$, which must be an integer. 

\noindent C. Since $s_{x,y}$ is determined by the vertex-types of $G$, it must remain the same in $\overline{G}$, that is, $s_{x,y} = \overline{s}_{x,y}$; similarly, $s_{\overline{x},\overline{y}} = \overline{s}_{\overline{x},\overline{y}}$. Now a vertex of degree $x$ (or $y$) in $G$ has degree $\overline{x}$ (or $\overline{y}$) in $\overline{G}$. So if there are $p$ unordered pairs of vertices $\{v,w\}$ with $d(v) = x$ and $d(w) = y$, then $p - s_{x,y} = \overline{s}_{\overline{x},\overline{y}} =  s_{\overline{x},\overline{y}}$. Since $p = {r_y \choose 2}$ when $x = y$, and $p = r_x r_y$ otherwise, the result follows.

\noindent D. By C, $\frac{\scriptstyle 1}{\scriptstyle 2} r_d^{\,2}$ must be an integer, so $r_d$ is even for $d \not= \overline{d}$. When there are vertices of degree $(n-1)/2$, by B we must have $n = 4k+1$ for some $k$. By the first part, there are $2 r$ vertices with degree $d < (n-1)/2$, and therefore $2 r$ vertices with degree $d > (n-1)/2$, leaving $(4k+1)-2(2r) = 4(k-r) + 1$ vertices of degree $(n-1)/2$.

\noindent E. If $G$ is vertex oblique, the vertices of degree $d$ cannot all be adjacent to vertices of degree $\overline{d}$ only. So 
$\frac{\scriptstyle 1}{\scriptstyle 2} r_d^{\,2} = s_{d,\overline{d}} < dr_d$.

\noindent F. Clearly $G$ cannot be regular, so there must be at least two different degrees. If $n$ is odd, the number of different degrees must be odd, by A. Suppose $n$ is even and there are exactly two degrees, say $d < (n-1)/2$ and $\overline{d} > (n-1)/2$. Since $G$ is vertex-oblique, each vertex $v$ of degree $d$ must be adjacent to a distinct number $n_v$ of vertices of degree $\overline{d}$. Note that $0 \leq n_v \leq d$. If  $n_v = 0$ for some $v$, then in the complement $v$ would be a vertex of degree $\overline{d}$ that is adjacent to all vertices of degree $d$, so $\overline{G}$ would have no vertex $w$ of degree $d$ with $n_w = 0$. Similar reasoning excludes the case $n_v = d$, so we have $0 < n_v < d$ for every $v$. This means that there are less than $d < n/2 = r_d$ possible values of $n_v$, a contradiction.
\risnumend
\newline

By Lemma~\ref{basic-lem}.B, the smallest possible orders $n>1$ for a dually vertex-oblique graph are $4$ and $5$; but by part E these could not have any vertices of degree $d < (n-1)/2$, so no such graphs exist when $n<8$. We now construct graphs for every $n \equiv 0$ or $1 \pmod 4$, $n \geq 8$.

\section{Construction on $4k$ vertices}

\begin{figure}[h]
\begin{center}
\input{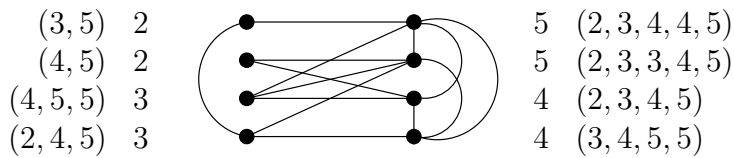}
\caption{A dually vertex-oblique graph on $8$ vertices.}
\label{fig-dvo-8}
\end{center}
\end{figure}

A dually vertex-oblique graph on $8$ vertices is shown in Figure~\ref{fig-dvo-8}, with the degree and vertex-type displayed next to each vertex. To verify this, one has to check that for every vertex with vertex-type $(x_1, \ldots, x_r)$, there is another vertex that is non-adjacent to vertices of degrees $(7-x_1, \ldots, 7-x_r)$.

\begin{figure}
\begin{center}
\input{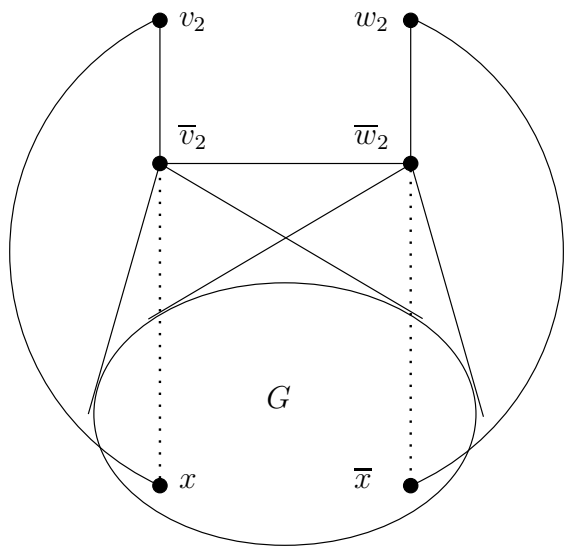}
\caption{Larger dually vertex-oblique graphs from smaller ones.}
\label{fig-dvo-n}
\end{center}
\end{figure}

Given a dually vertex-oblique graph $G$ on $n = 4k$ vertices, we now show how to construct $G'$ on $n+4 = 4k'$ vertices, where $k' := k+1$ (see Figure~\ref{fig-dvo-n}). We add vertices $v_2, w_2$, that will have degree $2$, and  $\overline{v}_2, \overline{w}_2$ that will have degree $2k' - 2 = 2k+2$. Moreover, the new vertices induce a $P_4$, in a manner reminiscent of Akiyama and Harary's~\cite{ag} method of producing larger self-complementary graphs.

We pick an arbitrary vertex $x \in V(G)$, and let $\overline{x}$ be the (unique) vertex such that $\overline{t}(\overline{x}) = t(x)$. Note that if $x$ has degree $d$, then $\overline{x}$ has degree $4k-1-d \not= d$. We make $v_2$ adjacent to $\overline{v}_2$ and $x$, $w_2$ adjacent to $\overline{w}_2$ and $\overline{x}$. Meanwhile, $\overline{v}_2$ is adjacent to $v_2, \overline{w}_2$ and $V(G) \setminus x$; and $\overline{w}_2$ is adjacent to $w_2, \overline{v}_2$ and $V(G) \setminus \overline{x}$.

A vertex of degree $d$ now has degree $d' := d+2$; this also means that a vertex of the complementary degree $\overline{d} := 4k-1-d$ now still has complementary degree $\overline{d'} = 4k'-1-d = \overline{d}+2$. The degrees in $V(G)$ now range between at least $4$ and at most $4k$ (by Lemma~\ref{basic-lem}.E), and thus the degrees (and vertex-types) in $V(G)$ are distinct from those of the new vertices. 

If $u \not\in \{x, \overline{x}\}$ had vertex-type $t(u) = (d_1, \ldots, d_r)$, in $G'$ it has type $t'(u) = (d_1', \ldots, d_r', 4k'-2, 4k'-2)$. The unique vertex $\overline{u}$ such that $\overline{t}(\overline{u}) = t(u)$ was non-adjacent in $G$ to vertices of degrees $\overline{d_1}, \ldots, \overline{d_r}$; so in $G'$ it is non-adjacent to vertices of degrees $\overline{d'_1}, \ldots, \overline{d'_r}$, as well as two vertices of degree $2 = \overline{4k'-2}$.

Meanwhile $t(x) = \overline{t}(\overline{x}) = (f_1, \ldots, f_d)$ becomes $t'(x) = \overline{t'}(\overline{x}) = (2, f_1+2, \ldots, f_d+2, 4k'-2)$. Distinct vertex-types in $G$ therefore result in distinct vertex-types in $G'$, and complementary vertex-types result in complementary vertex-types.
Moreover, $t(v_2) = \overline{t}(\overline{v}_2) = (2, d)$, and $t(w_2) = \overline{t}(\overline{w}_2) = (2, \overline{d})$, so the types of the new vertices are also distinct and complementary, and thus $G$ is dually vertex-oblique.
\newline

A graph $G$ is {\em split} if its vertex set partitions into $L \cup R$ (the ``left'' and ``right'' vertices), where $G[L]$ is edgeless and $G[R]$ is complete. Our constructions are close to being split graphs, with the vertex-set partitioning into vertices of degree less than $2k$ and vertices of degree at least $2k$. With a little more effort we can construct examples (for $n \geq 12$) that are actually split graphs. 

We will construct an appropriate bipartite graph with partition $L \cup R$ and show that, if we add edges to make $R$ induce a clique, the resulting graph is dually vertex-oblique. If $B$ is bipartite, with partition $L \cup R$, its {\em bipartite complement} is the graph $\tilde{B}$ with $V(\tilde{B}) := V(B)$, $E(\tilde{B}) := \{uv \mid u \in L, v \in R, uv \not\in E(B)\}$. The vertex-type of $v$ in $\tilde{B}$ is $\tilde{t}(v)$. A {\em dually semi-vertex-oblique} graph is a bipartite graph $B$ with $L = \{\ell_1, \ldots, \ell_{2k}\}, R = \{r_1, \ldots, r_{2k}\}$, such that:
\begin{itemize}
\item[(i)] $\{t(\ell_1), \ldots, t(\ell_{2k})\}$ contains no repetitions
\item[(ii)] $\{t(\ell_1), \ldots, t(\ell_{2k})\} = \{t(r_1), \ldots, t(r_{2k})\}$
\item[(iii)] $\{\tilde{t}(\ell_1), \ldots, \tilde{t}(\ell_{2k})\} = \{t(\ell_1), \ldots, t(\ell_{2k})\}$, and (thus) \\
$\{\tilde{t}(r_1), \ldots, \tilde{t}(r_{2k})\} = \{t(r_1), \ldots, t(r_{2k})\}$.
\end{itemize}

If $\ell_1$, say, had degree $2k$, then in $\tilde{B}$ it would have degree $0$, so by conditions (ii) and (iii) there must be a vertex of degree $0$ in $R$, a contradiction. So the minimum degree in $L$ is at least $1$, and by (iii) the maximum degree is at most $2k-1$, and similarly for $R$.

We now add edges to make $R$ induce a clique, giving us a split graph $G$. The degree of any vertex $r_j$ jumps up by $2k-1$, so its degree becomes at least $2k$; thus the degrees (and vertex-types) of vertices in $R$ become distinct from those in $L$. If  $t(r_i)$ differed from $t(r_j)$ in the number of entries equal to $d$, where $d < 2k$, then in $G$ both types get $2k-1$ new entries that are all at least $2k$, but they still differ in the number of entries equal to $d$. If $t(\ell_i)$ differed from $t(\ell_j)$ in the number of entries equal to $d$, in $G$ they will differ in the number of entries equal to $d+2k-1$. Thus $G$ is vertex-oblique, and from (iii) we can see that $\overline{G}$ has the same vertex-types as $G$.

\begin{figure}
\begin{center}
\input{dually_2.pstex_t}
\caption{Dually semi-vertex-oblique graphs on $12$ and $16$ vertices.}
\label{fig-sdvo-12-16}
\end{center}
\end{figure}

Dually semi-vertex-oblique graphs on $12$ and $16$ vertices are shown in Figure~\ref{fig-sdvo-12-16}, with their degrees and vertex-types.
To verify condition (iii), one has to check that for every vertex on the left with vertex-type $(x_1, \ldots, x_r)$, there is another vertex that is non-adjacent to vertices (on the right) of degrees $(2k-x_1, \ldots, 2k-x_r)$.

\begin{figure}
\begin{center}
\input{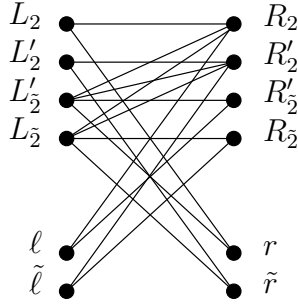}
\caption{Adding vertices to make larger semi-vertex-oblique graphs.}
\label{fig-sdvo-n}
\end{center}
\end{figure}

Given a dually semi-vertex-oblique graph $B$ on $n = 4k$ vertices, we now show how to construct $B'$ on $n+8 = 4k'$ vertices, where $k' := k+2$. We add vertices $L_2, L'_2, L_{\tilde{2}}, L'_{\tilde{2}}$ on the left, and $R_2, R'_2, R_{\tilde{2}}, R'_{\tilde{2}}$ on the right. The vertices with subscript $2$ will have degree $2$, those with subscript $\tilde{2}$ will have degree $2k' - 2 = 2k+2$. See Figure~\ref{fig-sdvo-n} for a sketch of the new vertices and their adjacencies to each other and to the vertices $\ell, \tilde{\ell}, r, \tilde{r}$ (described below); the vertices with subscript $\tilde{2}$ are also adjacent to all other vertices on the opposite side.

By (iii) we can find two vertices $r, \tilde{r}$ such that $\tilde{t}(\tilde{r}) = t(r)$; in particular, if $r$ has degree $d$, then $\tilde{r}$ has degree $2k-d$. We make $L_2$ adjacent to $R_{\tilde{2}}$ and $r$, $L_{\tilde{2}}$ adjacent to $R_2, R_{\tilde{2}}, R'_{\tilde{2}}$ and $\{r_i \not= r\}$. Similarly $L'_2$ is adjacent to $R'_{\tilde{2}}$ and $\tilde{r}$, $L'_{\tilde{2}}$ adjacent to $R'_2, R_{\tilde{2}}, R'_{\tilde{2}}$ and $\{r_i \not= \tilde{r}\}$. 

By (ii) there are (unique) vertices $\ell, \tilde{\ell}$ with $t(\ell) = t(r), t(\tilde{\ell}) = t(\tilde{r})$. The adjacencies for $R_2, R'_2, R_{\tilde{2}}, R'_{\tilde{2}}$ are defined as above: 
$N(R_2) := \{L_{\tilde{2}}, \tilde{\ell}\}$, $N(R_{\tilde{2}}) := \{L_2, L_{\tilde{2}}, L'_{\tilde{2}}\} \cup \{\ell_i \not= \tilde{\ell}\}$, $N(R'_2) := \{L'_{\tilde{2}}, \ell\}$, $N(R'_{\tilde{2}}) := \{L'_2, L_{\tilde{2}}, L'_{\tilde{2}}\} \cup \{\ell_i \not= \ell \}$.

The adjacencies of the new vertices are well-defined, and the construction is symmetric (as far as degrees and vertex-types go) with respect to $L$ and $R$, so (ii) holds. The degrees of every $\ell_i$ increase by $2$ (so they now range between at least $3$ and at most $2k+1$). If $t(\ell_i)$ differed from $t(\ell_j)$ in the number of entries equal to $d$, in $B'$ they will differ in the number of entries equal to $d+2$. By construction, $t(L_2), t(L'_2), t(L_{\tilde{2}})$ and $t(L'_{\tilde{2}}$ are distinct from each other (and from the $t(\ell_i)$'s, because of their degrees). Thus (i) holds.

In the bipartite complement, $N(L_2) = \{R_2, R'_2, R'_{\tilde{2}}\} \cup \{r_i \not= r\}$, $N(L_{\tilde{2}}) = \{R'_2, r\}$, $N(L'_2) = \{R_2, R'_2, R_{\tilde{2}}\} \cup \{r_i \not= \tilde{r}\}$ and $N(L'_{\tilde{2}}) = \{R_2, \tilde{r}\}$. The neighbourhoods of $R_2, R'_2, R_{\tilde{2}}$ and $R'_{\tilde{2}}$ are changed similarly. Recall also that $\tilde{t}(\tilde{\ell}) = t(\ell)$ and $\tilde{t}(\tilde{r}) = t(r)$. Thus $\tilde{B'}$ is obtained from $\tilde{B}$ in the same way as we obtained $B'$ from $B$ (with the roles of $L_2$ and $L'_{\tilde{2}}$ interchanged, and similarly for $L'_2$ and $L_{\tilde{2}}$, $\ell$ and $\tilde{\ell}$, and so on). Since $B$ satisfied (iii), $B'$ does too.

\section{Construction on $4k+1$ vertices}
Take a dually vertex-oblique graph $G$ on $4k$ vertices, and introduce a new vertex $u_0$ that is adjacent to the $2k$ vertices with degree $d \geq 2k$. We claim that the  resulting graph $G'$ of order $n' = 4k+1$ is again dually vertex-oblique. Note that if $G$ was a split graph, then the high-degree vertices must have formed a clique, and thus $G'$ will also be split.

A vertex of degree $d := d_G(v)$ in $G$ has degree $d' := d_{G'}(v)$ in $G'$. So $d' = d$ if $d < 2k$, and $d' = d+1$ if $d \geq 2k$. 
If $v$ and $w$ had complementary degrees in $G$, that is, $d_G(v) + d_G(w) = n-1 = 4k-1$, then $d_{G'}(v) + d_{G'}(w) = n'-1 = 4k$, so $v$ and $w$ still have complementary degrees in $G'$; this means that if $\overline{d} = f$ in $G$, then $\overline{d'} = f'$ in $G'$. Also, $d_{G'}(u_0) + d_{G'}(u_0) = 2k+2k = n'-1$.

In $G'$, $u_0$ will be the unique vertex of degree $2k$; and if $v,w,$ were adjacent in $G$ to different numbers of vertices of degree $d$, then in $G'$ they are adjacent to different numbers of vertices of degree $d'$; so $G'$ is vertex-oblique.

Since $G$ is dually vertex-oblique, for every vertex $v$ there is a unique vertex $\overline{v}$ with $\overline{t}(\overline{v}) = t(v)$. If $t(v) = (x_1, \ldots, x_r)$, then $\overline{v}$ must have non-neighbours in $G$ of degrees $\overline{x_1}, \ldots, \overline{x_r}$. If $r < 2k$, then $u_0$ is adjacent to $\overline{v}$ but not to $v$, so in $G'$ $v$ has vertex-type $t'(v) = (x'_1, \ldots, x'_r)$, and $\overline{v}$ has non-neighbours of degrees $\overline{x_1'}, \ldots, \overline{x_r'}$; thus $\overline{t}'(\overline{v}) = t'(v)$. If $r > 2k$ then $u_0$ is adjacent to $v$ in $G'$, and to $\overline{v}$ in $\overline{G'}$; thus $t'(v) = (x'_1, \ldots, 2k, \ldots, x'_r) = \overline{t}'(\overline{v})$. Finally $t'(u_0) = \overline{t}'(u_0)$, so $G'$ is dually vertex-oblique.

\section{Vertex-type sequences: \\ uniqueness and non-uniqueness}

The {\em degree sequence} of a graph on $n$ vertices is the sequence $d_1 \geq \cdots \geq d_n$ of its degrees (see footnote~\ref{foot-1}, p.~\pageref{foot-1}). The {\em vertex-type sequence} is the sequence $t_1 \succeq \cdots \succeq t_n$ of vertex-types, where $t_i \succ t_j$ if $t_i$ is longer than $t_j$, or if $t_i$ and $t_j$ have the same length and $t_i$ is lexicographically larger than $t_j$. 
$G_d$ is the subgraph of $G$ induced by vertices of degree $d$, and (for $p \not= q$) $G_{p,q}$ is the bipartite subgraph induced by edges joining a vertex of degree $p$ to a vertex of degree $q$. 

Some graphs, such as complete graphs, edgeless graphs and matchings, have unique degree sequence (that is, no other graph has the same degree sequence) and, thus, unique vertex-type sequence. If $G$ is dually vertex-oblique, then by definition its complement shares the same vertex-type sequence, and is not isomorphic to $G$ because self-complementary graphs have non-trivial automorphisms. But could this complementary pair be the unique graphs with that vertex-type sequence? We show here that the answer is always `No', but that there are infinitely many super vertex-oblique graphs with unique vertex-type sequence.
\newline

The key to the proofs is a restricted switching operation. A {\em switch} is the replacement of edges $v_0 w_0, v_1 w_1$, with new edges $v_0 w_1, v_1 w_0$ (that is, $v_0 w_1$ and $v_1 w_0$ did not appear in the original graph); this does not change the degree of any vertex, but may change the vertex-types. A {\em $(d,d')$-switch} (or just `restricted switch', when $d$ and $d'$ are not specified) is a switch where $v_0$ and $v_1$ both have degree $d$, and $w_0$ and $w_1$ both have degree $d'$ (possibly equal to $d$); such a switch does not change the type of any vertex. In a bipartite graph, a switch {\em respects the bipartition} if $v_0, v_1$, are in the same part, and (thus) $w_0, w_1$, are in the opposite part.

\resnum{Theorem}
{For any dually vertex-oblique graph $G$, there is a graph $H \not\in \{G, \overline{G}\}$ with the same vertex-type sequence as $G$.
}
We will establish:
\newline

\noindent{\sc Claim.} For any degree $d \not= (n-1)/2$, $G$ has distinct vertices $v_0$, $v_1$, of degree $d$, and $w_0, w_1$, of degree $\overline{d}$, such that $v_0 w_0, v_1 w_1 \in E(G)$, $v_0 w_1, v_1 w_0 \not\in E(G)$. 
\newline

The result follows from the claim since we can then perform a $(d, \overline{d})$-switch which gives us another graph $H$ without changing the type of any vertex. 
Since $G$ has trivial automorphism group, $H \not\cong G$, and we will show that $H \not\cong \overline{G}$.

Let $x$ be any vertex not in $\{v_0, v_1, w_0, w_1\}$. Let $\overline{x}$ be the unique vertex that has the same vertex-type in $\overline{G}$ as $x$ has in $G$. If $\overline{x} \not= x$ (possibly $\overline{x} \in \{v_0, v_1, w_0, w_1\}$), note that $x$ and $\overline{x}$ are adjacent in $H$ iff they are adjacent in $G$ iff they are \emph{not} adjacent in $\overline{G}$. If $x = \overline{x}$, take another vertex $y \not\in \{x, v_0, v_1, w_0, w_1\}$; note that $x$ and $y$ exist by the remark after Lemma~\ref{basic-lem}. If $y \not= \overline{y}$ we are done, otherwise note that $x$ and $y$ are adjacent in $H$ iff they are adjacent in $G$ iff $\overline{x}=x$ and $\overline{y}=y$ are not adjacent in $\overline{G}$.
\newline

We now turn to proving the Claim, which is equivalent to saying that the be the bipartite graph $G_{d, \overline{d}}$ has an induced $2 K_2$. If there is any vertex $z$ of degree $d$ that is adjacent to no (or all) vertices of degree $\overline{d}$, then in $\overline{G}$ $z$ would be a vertex of degree $\overline{d}$ adjacent to all (or no) vertices of degree $d$, contradicting the fact that $G$ and $\overline{G}$ have the same vertex-type sequence. So in $G_{d, \overline{d}}$ every vertex has at least one neighbour and one non-neighbour from the opposite part.

Let $G_{d, \overline{d}}$ have bipartition $D \cup \overline{D}$. In what follows, $x_i$ will be a vertex in $D$, $N_i \subseteq \overline{D}$ the set of its neighbours, and $N'_i := \overline{D} \setminus N_i$; $N_i \not= \emptyset \not= N'_i$ by the previous argument. Take an arbitrary vertex $x_0 \in D$. Pick a vertex $y_1 \in N'_0$, and let $y_1$ be adjacent to some vertex $x_1$; clearly $x_1 \not= x_0$. If there is a vertex $\tilde{y}_0 \in N_0$ such that $\{x_0,x_1, \tilde{y}_0, y_1\}$ induce a $2K_2$, we are done; otherwise, $x_1$ is adjacent to all of $N_0$, as well as $y_1$, so $N_1 \supsetneq N_1$, and $N'_1 \subsetneq N'_0$. Pick a vertex $y_2 \in N_1$, and let $y_2$ be adjacent to $x_2$; as before, $x_2 \not= x_1$, and either there is $\tilde{y}_1 \in N_1$ such that $\{x_1,x_2, \tilde{y}_1, y_2\}$ induce a $2K_2$, or $N'_2 \subsetneq N'_1$. Repeating this procedure we must eventually find an induced $2 K_2$, since $N'_i$ can never be empty.
\resnumend
\newline

A graph $G$ can be transformed by switches into any other graph $H$ with the same degree sequence. If $H$ even has the same vertex-type sequence as $G$, then we will show how to achieve this using only restricted switches.

Suppose the vertices of a graph $G$ are labeled $v_1, \ldots, v_n$, with $\Delta = d(v_1) \geq \cdots \geq d(v_n)$. By switching, we can transform $G$ into a canonical labeled graph $G_0$ that is determined completely by the degree sequence (the first step in this recursive process is to use switches to make $v_1$ adjacent to $v_2, \ldots, v_{\Delta+1}$); any other labeled graph $H$ with the same vertex-set and the same degree sequence (i.e. $d_G(v_i) = d_H(v_i)$ for all $i$) can also be transformed into $G_0$. These ideas, and analogous ones for bipartite graphs, give us:

\stmt{Theorem~\cite{egg, fulk}}
{If $G$, $H$, are two labeled graphs with the same degree sequence, then $G$ can be obtained from $H$ by a sequence of switches. Moroever, if $G$ and $H$ are bipartite, the switches respect the bipartition.
}
\stmtend
\newline

We use this to prove the next result, that has probably also appeared in~\cite{zver}:

\resnum{Theorem}
{If $G$, $H$, are two labeled graphs with the same vertex-type sequence, then $G$ can be obtained from $H$ by a sequence of restricted switches.
}
The vertex-type sequence clearly determines the degree sequence. Moreover, for every degree $d$, the subgraphs $G_d$ and $H_d$ have the same vertex-set and the same degree sequence, since this is also determined by the vertex-types; we can therefore transform $G_d$ into $H_d$ by a sequence of switches; note that in $G$ these are just $(d,d)$-switches. Similarly, for every $p \not=q$ in the degree sequence, the bipartite graphs $G_{p,q}$ and $H_{p,q}$ have the same vertex-set, the same bipartition, and the same degrees, so we can transform $G_{p,q}$ into $H_{p,q}$ by switches; moreover, we can use switches that respect the bipartition, and these will be valid $(p,q)$-switches in $G$ even though $G_{p,q}$ is not a vertex-induced subgraph of $G$.
\resnumend

\stmtnum{Corollary}
{If no $G_d$ and no $G_{p,q}$ contains an induced $2 K_2$, then $G$ has unique vertex-type sequence. In particular, if every degree appears at most once in $G$, except for some degree that appears at most three times, then $G$ has unique vertex-type sequence.
}
\stmtnumend
\newline

The converse of the corollary is not true (the matchings are a counterexample), because a restricted switch may give us a graph $G'$ isomorphic to $G$. But it can be used to show, for example, that  the super vertex-oblique graphs $G_1^6, G_1^7, G_2^8, G_2^9$, in~\cite{sch} have unique vertex-type sequence. In particular, every degree appears exactly once in $G_2^8$, except for five vertices of the same degree that induce a graph with only one edge; it can be checked that applying Construction 1 of~\cite{sch} with $k = 1$ preserves these properties, and it is shown in that paper that the result is again connected and super vertex-oblique. We thus have:

\stmtnum{Corollary}
{There are infinitely many connected super vertex-oblique graphs with a unique vertex-type sequence.
}
\stmtnumend

\section{Recognising degree and vertex-type sequences\label{sec-deg-seq}}
A graph $G$ {\em realises} its degree sequence, and its vertex-type sequence.
Erd\H{o}s and Gallai%
\footnote{Several authors have given different characterisations of degree sequences of graphs.}~\cite{em1}
showed that a sequence $d_1 \geq \cdots \geq d_n$ is realised by some graph if and only if, for $r = 1,\ldots,n-1$, we have
\[ \sum_{i=1}^{r} d_i \leq r(r-1) + \sum_{j=r+1}^{n} \min(r, d_j). \]

If $G$ is a graph with the same degree sequence 
as $\overline{G}$, and $r_d$ is the number of vertices of degree $d$, then:
\begin{itemize}
\item[(*)] $d_i + d_{n-i+1} = n-1$, for $i = 1, \ldots, n$;
\item[(**)] $r_d$ is even for all $d$, except for $r_{(n-1)/2} \equiv 1 \pmod 4$.
\end{itemize}
Clapham and Kleitman~\cite{do} showed by construction that every sequence that satisfies (*), (**) and the Erd\H{o}s-Gallai conditions, is realised by a self-complementary graph. However (Lemma~\ref{basic-lem}.F), not all such sequences are realised by a dually vertex-oblique graph. It would be interesting to characterise the degree sequences of dually vertex-oblique graphs. One might also ask similar questions about vertex-type sequences:

\stmtnobox{Problem}
{A. When is a sequence (of sequences of positive integers) the vertex-type sequence of some graph?

\noindent B. Characterise the degree sequences and vertex-type sequences of:
\begin{itemize}
\item vertex-oblique graphs,
\item super vertex-oblique graphs, and
\item dually vertex-oblique graphs.
\end{itemize}
}
\stmtnoboxend

The Erd\H{o}s-Gallai results on degree-sequences, together with the Gale-Ryser conditions explained below, lead to an efficient algorithm to solve the vertex-type sequence problems; if the sequence is realised by some graph, the algorithm can also be made to construct an example. However, we would like a more succint characterisation similar to that of Erd\H{o}s-Gallai, Clapham-Kleitman or Gale-Ryser, especially as this might shed light on the degree sequence problems.

Gale~\cite{gale} and Ryser~\cite{ryser} showed that sequences $p_1 \geq \cdots \geq p_m$ and $q_1 \geq \cdots \geq q_n$ are the degrees of a bipartite graph $B$ (with the $p_i$'s being degrees on one side, and the $q_j$'s the degrees on the other side) if and only if, for $r = 1,\ldots,n-1$:
\[ \sum_{i=1}^m \min\{r, p_i\} \geq \sum_{j=1}^r q_j .\]

Given the vertex-type sequence of a graph $G$, we can recover the degree sequence, and compute the vertex-types of $\overline{G}$ (as noted at the beginning of the introduction); it is then straightforward to check whether $G$ is (super or dually) vertex-oblique. So we turn our attention to Problem A.

If we want to check whether a given sequence is actually the vertex-type sequence of some graph $G$, we recover the degree-sequences of the $G_d$'s and $G_{p,q}$'s (for all $d,p,q,$ in the degree-sequence of $G$), and check the Erd\H{o}s-Gallai and Gale-Ryser conditions, respectively. If the conditions are not all satisfied, we have a contradiction; otherwise, we can construct $G_d$'s and $G_{p,q}$'s that together give us a graph with the given vertex-type sequence.

\section{Other open problems}
In a self-complementary graph of order $4k+1$, one can always remove an appropriate vertex to get a  self-complementary graph of order $4k$. It is not clear whether an analogous claim is true for dually vertex-oblique graphs.

\stmtnobox{Problem}
{Is there a dually vertex oblique graph on $4k+1$ vertices, such that removing any vertex of degree $2k$ leaves a subgraph $H$ such that (a) $H$ does not have the same vertex-types as its complement, or (b) $H$ is not vertex-oblique, or both (a) and (b)?
}
\stmtnoboxend
\newline

For any fixed $k$, Schreyer et al.~\cite{sch} constructed super vertex-oblique graphs that were $k$-connected, with $k$-connected complements. Our examples of dually vertex-oblique graphs have vertices of degree $2$, and thus connectivity at most $2$. 

\stmtnobox{Problem}
{Are there (complementary pairs of) dually vertex-oblique graphs of arbitrarily high connectivity?
}
\stmtnoboxend

\section{Acknowledgements}
My studies in Canada are fully funded by the
Canadian government through a Canadian Commonwealth Scholarship.

I would like to thank Jens Schreyer and Hansjoachim Walther for introducing me to the concept of vertex-oblique graphs at the Cycles and Colourings workshop in Star\'{a} Lesn\'{a} (Slovakia, 2002), and the Graph Theory conference at Czorstyn (Poland, 2002). Participation in these conferences was funded by the University of Waterloo.

\end{document}